\pdfoutput=1
\RequirePackage{ifpdf}
\ifpdf 
\documentclass[pdftex]{sigma}
\else
\documentclass{sigma}
\fi

\newtheorem{criterion}{Criterion}
\def\Poch#1#2#3{{\left(#1;#2\right)\sb{#3}}}

\begin{document}

\allowdisplaybreaks

\renewcommand{\thefootnote}{$\star$}

\renewcommand{\PaperNumber}{015}

\FirstPageHeading

\ShortArticleName{On the Smoothness of the Noncommutative Pillow and Quantum Teardrops}

\ArticleName{On the Smoothness of the Noncommutative Pillow\\
and Quantum Teardrops\footnote{This paper is a~contribution to the Special Issue on Noncommutative Geometry
and Quantum Groups in honor of Marc A.~Rief\/fel.
The full collection is available at
\href{http://www.emis.de/journals/SIGMA/Rieffel.html}{http://www.emis.de/journals/SIGMA/Rieffel.html}}}

\Author{Tomasz BRZEZI\'NSKI}

\AuthorNameForHeading{T.~Brzezi\'nski}

\Address{Department of Mathematics, Swansea University, Singleton Park, Swansea SA2 8PP, UK}
\Email{\href{mailto:T.Brzezinski@swansea.ac.uk}{T.Brzezinski@swansea.ac.uk}}

\ArticleDates{Received December 03, 2013, in f\/inal form February 09, 2014; Published online February 14, 2014}

\Abstract{Recent results by Kr\"ahmer~[\textit{Israel~J. Math.} \textbf{189} (2012), 237--266] on smoothness of Hopf--Galois extensions and by
Liu~[arXiv:1304.7117] on smoothness of generalized Weyl algebras are used to prove that the coordinate
algebras of the noncommutative pillow orbifold~[\textit{Internat.~J. Math.} \textbf{2} (1991), 139--166], quantum teardrops ${\mathcal O}({\mathbb
W}{\mathbb P}_q(1,l))$~[\textit{Comm. Math. Phys.} \textbf{316} (2012), 151--170], quantum lens spaces ${\mathcal O}(L_q(l;1,l))$~[\textit{Pacific~J. Math.} \textbf{211} (2003), 249--263],
the quantum Seifert manifold ${\mathcal O}(\Sigma_q^3)$~[\textit{J.~Geom. Phys.} \textbf{62} (2012), 1097--1107], quantum real weighted projective
planes ${\mathcal O}({\mathbb R}{\mathbb P}_q^2(l;\pm))$~[\textit{PoS Proc. Sci.}  (2012), PoS(CORFU2011), 055, 10~pages] and quantum Seifert lens spaces
${\mathcal O}(\Sigma_q^3(l;-))$~[\textit{Axioms} \textbf{1} (2012), 201--225] are homologically smooth in the sense that as their own
bimodules they admit f\/initely generated projective resolutions of f\/inite length.}


\Classification{58B32; 58B34}

\begin{flushright}
\em Dedicated to Marc Rieffel on the occasion of his 75th birthday.
\end{flushright}

\renewcommand{\thefootnote}{\arabic{footnote}}
\setcounter{footnote}{0}

\section{Introduction and tools}

This note is intended to illustrate the claim that (often)
a $q$-deformation of a~non-smooth classical variety or an orbifold produces an algebra which
has properties of the coordinate algebra of a~non-commutative smooth variety or manifold.
More precisely, we say that an algebra $B$ (over an algebraically closed f\/ield ${\mathbb K}$) is {\em
homologically smooth} or simply {\em smooth} provided that as a~$B$-bimodule it has a~f\/initely generated
projective resolution of f\/inite length; see~\cite[Erratum]{Van:rel}.
We prove that several classes of examples of coordinate algebras of $q$-deformed orbifolds are
homologically smooth.
To achieve this aim we use techniques developed in~\cite{Kra:Hoc}, which are applicable to principal
comodule algebras~\cite{BrzHaj:Che}, and those of~\cite{Liu:hom}, which are applicable to generalized Weyl
algebras~\cite{Bav:gweyl}.
We summarize these presently.

In~\cite[Corollary~6]{Kra:Hoc} Kr\"ahmer gives a~criterion for smoothness of quantum homogeneous spaces,
which through an immediate extension to a~more general class of Hopf--Galois extensions and then
specif\/ication to strongly group-graded algebras, provides us with a~tool for showing the smoothness of
the noncommutative pillow algebra studied in~\cite{BraEll:sph}, the quantum lens space algebras ${\mathcal
O}(L_q(l;1,l))$
introduced in~\cite{HonSzy:len}, the quantum teardrop algebras~\cite{BrzFai:tea},
the coordinate algebras of quantum real weighted
projective planes ${\mathcal O}({\mathbb R}{\mathbb P}_q(l;-))$ def\/ined in~\cite{Brz:Sei}, the quantum
Seifert manifold ${\mathcal O}(\Sigma_q^3)$~\cite{BrzZie:pri} and the quantum Seifert lens spaces
${\mathcal O}(\Sigma_q^3(l;-))$~\cite{BrzFai:rps}.
More specif\/ically, let $G$ be a~group (with the neutral element~$e$).
A $G$-graded algebra $A=\oplus_{g\in G} A_g$ is said to be {\em strongly graded} if, for all $g,h\in G$,
$A_gA_h = A_{gh}$.
For such algebras, Kr\"ahmer's criterion for smoothness takes the following form.
\begin{criterion}
\label{cri.smooth}
Let $A$ be a~strongly $G$-graded algebra and set $B= A_e$.
If the enveloping algebra ${\mathcal E} A:= A\otimes A^{\rm op}$ of $A$ is left Noetherian of finite global
dimension, then the enveloping algebra of~$B$ is also left Noetherian of finite global dimension.
Consequently, $B$ is a~homologically smooth algebra.
\end{criterion}

Although~\cite[Corollary~6]{Kra:Hoc} is formulated for quantum homogeneous spaces obtained via surjective
homomorphism of Hopf algebras with a~cosemisimple codomain, the proof extends immediately to all faithfully
f\/lat Hopf--Galois extensions $B\subseteq A$ or principal comodule algebras such that $B$ is a~direct
summand of $A$ as a~$B$-bimodule.
In the case of $G$-graded algebras, $B=A_e$ is a~direct summand of $A$ as a~$B$-bimodule,
and~\cite[Proposition~AI.3.6]{NasVan:gra} ensures that a~strongly $G$-graded algebra is a~principal
${\mathbb K} G$-comodule algebra.
In particular $A$ is projective and faithfully f\/lat as a~left and right $B$-module.
The way Criterion~\ref{cri.smooth} is stated indicates its iterative nature which is implicit in the proof
of~\cite[Corollary~6]{Kra:Hoc}.
Note in passing that the assumption about the dimension of ${\mathcal E} A$ is not necessary to conclude
that ${\mathcal E} B$ is left Noetherian.

An ef\/fective way of checking whether a~graded $G$-algebra $A$ is strongly graded is described
in~\cite[Section~AI.3.2]{NasVan:gra}:
\begin{lemma}\label{lem.strong}
$A= \oplus_{g\in G} A_g$ is strongly graded if and only if there exists a~function $\omega: G \to A\otimes A$ such that

$(a)$ for all $g\in G$, $\omega(g) \in A_{g^{-1}}\otimes A_g$,

$(b)$ for all $g\in G$, $\mu\circ\omega(g) = 1$, where $\mu$ is the multiplication map of $A$.
\end{lemma}

Furthermore, if $G$ is a~cyclic group, then conditions $(a)$ and $(b)$ need only be checked for a~generator $g$ of~$G$.
If $\omega(g)=\sum\limits_i\omega'_i\otimes \omega''_i$, satisf\/ies $(a)$ and $(b)$, then $\omega$ is def\/ined by
setting $\omega(e) =1\otimes 1$~and
\begin{gather*}
\omega\big(g^{n+1}\big)=\sum_i\omega'_i\omega(g^n)\omega''_i,
\qquad
\text{for~all}
\quad
n>0.
\end{gather*}
A function $\omega$ satisfying conditions $(a)$ and $(b)$ in Lemma~\ref{lem.strong} is a~predecessor of a~{\em
strong connection form} on a~principal comodule algebra;
see~\cite{BrzHaj:Che, DabGro:str,Haj:str}.

Let $R$ be an algebra, let $p$ be an element of the centre of $R$ and let $\pi$ be an automorphism of~$R$.
The {\em $($degree-one$)$ generalized Weyl algebra} $R(\pi,p)$ is the extension of $R$ by genera\-tors~$x_+$,~$x_-$
subject to the relations, for all $r\in R$,
\begin{gather*}
x_-x_+=p,
\qquad
x_+x_-=\pi(p),
\qquad
x_{\pm}r=\pi^{\pm1}(r)x_\pm;
\end{gather*}
see~\cite{Bav:gweyl}.
In~\cite[Theorem~4.5]{Liu:hom} Liu gives the following criterion of smoothness of a~generalized Weyl
algebra over the polynomial algebra.
\begin{criterion}
\label{cri.Weyl}
Let $R={\mathbb K}[a]$ be a~polynomial algebra and an automorphism $\pi: {\mathbb K}[a]\to {\mathbb K}[a]$
be determined by $\pi(a) = \kappa a+\chi$.
Then the generalized Weyl algebra $R(\pi,p)$ is homologically smooth with homological dimension~$2$ if and
only if the polynomial $p\in {\mathbb K}[a]$ has no multiple roots.
\end{criterion}

Furthermore, Liu proves that if the smoothness Criterion~\ref{cri.Weyl} is satisf\/ied, then $A=R(\pi,p)$
is a~{\em twisted Calabi--Yau algebra} of dimension $2$ with the Nakayama (twisting) automorphism $\nu: A\to
A$ given by $\nu(x_\pm) = \kappa^{\pm 1} x_\pm$ and $\nu(a) = a$.
This means that the Hochschild cohomology of $A$ with values in its enveloping algebra is trivial in all
degrees except degree~2, where it is equal to $A$ with the $A$-bimodule structure $a\cdot b \cdot a' = ab\nu(a')$.

The reader should observe that, except for some special cases, the algebras described herein\-after are not
smooth whenever the deformation parameters $\lambda$ or $q$ are equal to 1.
By noting this they will fully grasp the main message of this note, namely that deformation may (and quite
often does) result in smoothing classically singular objects.

\section{Results}

Throughout we work with associative complex $*$-algebras with identity.
We write ${\mathcal E} A$ for the enveloping algebra $A\otimes A^{\rm op}$ of $A$.
We often use the {\em $q$-Pochhammer symbol} which, for an indeterminate $x$ and a~complex number $q$, is
def\/ined as
\begin{gather*}
\Poch xqn:=\prod_{m=0}^{n-1}\big(1-q^{m}x\big).
\end{gather*}

\subsection{The noncommutative pillow}

Let $\lambda = e^{2\pi i\theta}$, where $\theta$ is an irrational number.
Recall that the coordinate $*$-algebra ${\mathcal O}(\mathbb{T}^2_\theta)$ of the {\em noncommutative
torus} is generated by unitaries $U,V$, such that $UV=\lambda VU$; see~\cite{Rie:alg}.
The involutive algebra automorphism given by
\begin{gather*}
\sigma: \ {\mathcal O}\big(\mathbb{T}^2_\theta\big)\to{\mathcal O}\big(\mathbb{T}^2_\theta\big),
\qquad
U\mapsto U^{*},
\qquad
V\mapsto V^{*},
\end{gather*}
makes ${\mathcal O}(\mathbb{T}^2_\theta)$ into a~${\mathbb Z}_2$-graded algebra.
The f\/ixed point (or degree-zero) subalgebra ${\mathcal O}(P_\theta)$ is generated by $U+U^*$ and $V+V^*$.
It has been introduced and studied from a~topological point of view in~\cite{BraEll:sph} (see
also~\cite[Section~3.7]{EvaKaw:sym}) as a~deformation of the coordinate algebra of the {\em pillow
orbifold}~\cite[Chapter~13]{Thu:geo} (an orbifold rather than manifold since, classically, the ${\mathbb
Z}_2$-action determined by the automorphism $\sigma$ is not free).
\begin{theorem}
\label{thm.main.pil}
{}${\mathcal O}(\mathbb{T}^2_\theta)$ is a~strongly ${\mathbb Z}_2$-graded algebra and the noncommutative
pillow algebra ${\mathcal O}(P_\theta)$ is homologically smooth.
\end{theorem}
\begin{proof}
Set
\begin{gather*}
\hat{x}=U-U^*,
\qquad
\hat{y}=V-V^*,
\qquad
\hat{z}=UV^*-U^*V.
\end{gather*}
Note that $\sigma(\hat{x}) = - \hat{x}$, $\sigma(\hat{y}) = - \hat{y}$ and $\sigma(\hat{z}) = - \hat{z}$,
i.e.\ all these are homogeneous elements of~${\mathcal O}(\mathbb{T}^2_\theta)$ with the ${\mathbb
Z}_2$-degree $1$.
A straightforward calculation af\/f\/irms that these elements satisfy the following relation
\begin{gather*}
\hat{x}^2+\hat{y}^2-\bar{\lambda}\hat{z}^2-\hat{x}z\hat{y}=2\big(\bar{\lambda}^2-1\big),
\end{gather*}
where $z= UV^* +U^*V\in {\mathcal O}(P_\theta)$.
Therefore, the mapping $\omega: {\mathbb Z}_2 \to {\mathcal O}(\mathbb{T}^2_\theta)\otimes {\mathcal
O}(\mathbb{T}^2_\theta)$, def\/ined as $\omega(0) =1\otimes 1$ and
\begin{gather*}
\omega(1)=\frac{1}{2\big(\bar{\lambda}^2-1\big)}\left(\hat{x}\otimes\hat{x}+\hat{y}\otimes\hat{y}
-\bar{\lambda}\hat{z}\otimes\hat{z}-\hat{x}z\otimes\hat{y}\right),
\end{gather*}
satisf\/ies conditions $(a)$ and $(b)$ in Lemma~\ref{lem.strong}, and ${\mathcal O}(\mathbb{T}^2_\theta)$ is
a~strongly ${\mathbb Z}_2$-graded algebra.

Both ${\mathcal O}(\mathbb{T}^2_\theta)$ and ${\mathcal E}{\mathcal O}(\mathbb{T}^2_\theta)$ can be
understood as iterated skew Laurent polynomial rings and hence they are left Noetherian
by~\cite[Theorem~1.4.5]{McCRob:Noe}.
Furthermore, the global dimension of the latter is less than or equal to 4
by~\cite[Theorem~7.5.3]{McCRob:Noe}.
Therefore, the noncommutative pillow algebra ${\mathcal O}(P_\theta)$ is homologically smooth by
Criterion~\ref{cri.smooth}.
\end{proof}

In short, Theorem~\ref{thm.main.pil} means that for the irrational $\theta$ (or, more generally, for any
real $\theta\in (0,1)\setminus\{\frac12\}$) the action of ${\mathbb Z}_2$ on the noncommutative torus is
free despite the fact that the corresponding action on the classical level is not free.
The set of f\/ixed points corresponds to a~manifold rather than an orbifold.

\subsection{Quantum teardrops and lens spaces}
\label{sec.lens}
Here we deal with three (classes of) complex $*$-algebras given in terms of generators and relations.

The coordinate algebra of the quantum three-sphere, ${\mathcal O}(S^{3}_q)$, is generated by $\alpha$ and
$\beta$ such that
\begin{gather}
\alpha\beta=q\beta\alpha,
\qquad
\alpha\beta^*=q\beta^*\alpha,
\qquad
\beta\beta^*=\beta^*\beta,
\nonumber
\\
\alpha\alpha^*=\alpha^*\alpha+\big(q^{-2}-1\big)\beta\beta^*,
\qquad
\alpha\alpha^*+\beta\beta^*=1,
\label{su}
\end{gather}
where $q\in(0,1)$; see~\cite{Wor:com}.
For any positive integer $l$, the coordinate algebra of the quantum lens space ${\mathcal O}(L_q(l;1,l))$
is a~$*$-algebra generated by $c$ and $d$ subject to the following relations:
\begin{gather*}
cd=q^{l}dc,
\qquad\!
cd^*=q^{l}d^*c,
\qquad\!
dd^*=d^*d,
\qquad\!
cc^*=\Poch{dd^*}{q^2}l,
\qquad\!
c^*c=\Poch{q^{-2}dd^*}{q^{-2}}l,
\end{gather*}
see~\cite{HonSzy:len}.
Finally, for a~positive integer $l$, the coordinate algebra of the quantum teardrop ${\mathcal
O}(\mathbb{WP}_q(1,l))$ is the $*$-algebra generated by $a$ and $b$ subject to the following relations
\begin{gather*}
a^*=a,
\qquad
ab=q^{-2l}ba,
\qquad
bb^*=q^{2l}a\Poch a{q^2}l,
\qquad
b^*b=a\Poch{q^{-2}a}{q^{-2}}l;
\end{gather*}
see~\cite{BrzFai:tea}.
These algebras form a~tower ${\mathcal O}(\mathbb{WP}_q(1,l))\hookrightarrow {\mathcal
O}(L_q(l;1,l))\hookrightarrow {\mathcal O}(S^{3}_q)$ with embeddings $a\mapsto dd^*$, $b\mapsto cd$ and
$c\mapsto \alpha^l$, $d\mapsto \beta$, respectively.
We thus can and will think of ${\mathcal O}(\mathbb{WP}_q(1,l))$ and ${\mathcal O}(L_q(l;1,l))$ as
subalgebras of ${\mathcal O}(S^{3}_q)$.
{}${\mathcal O}(S^{3}_q)$ is a~${\mathbb Z}_l$-graded algebra with grading given by $\deg(\alpha) =1$,
$\deg(\alpha^*) = l-1$, $\deg(\beta)=\deg(\beta^*) =0$, and the above embedding identif\/ies the
degree-zero part of ${\mathcal O}(S^{3}_q)$ with ${\mathcal O}(L_q(l;1,l))$.
By~\cite[Theorem~3.3]{BrzFai:tea}, the latter is a~strongly ${\mathbb Z}$-graded algebra with grading
provided by $\deg(c) = \deg(d^*) = 1$, $\deg(c^*) = \deg(d) =-1$ and with the degree-zero part isomorphic
to ${\mathcal O}(\mathbb{WP}_q(1,l))$.

That ${\mathcal O}(\mathbb{WP}_q(1,l))$ is homologically smooth can be argued as follows.
{}${\mathcal O}(S^{3}_q)$ is a~coordinate algebra of the quantum group ${\rm SU}(2)$ and thus ${\mathcal
E}{\mathcal O}(S^{3}_q)$ is left Noetherian and has a~f\/inite global dimension; see~\cite{GooZha:hom}.
Hence, if it were a~strongly ${\mathbb Z}_l$-graded algebra, then ${\mathcal E}{\mathcal O}(L_q(l;1,l))$
would be left Noetherian and would have a~f\/inite global dimension (so, in particular ${\mathcal
O}(L_q(l;1,l))$ would be homologically smooth) by Criterion~\ref{cri.smooth}.
Since, in turn ${\mathcal O}(L_q(l;1,l))$ is a~strongly graded algebra, Criterion~\ref{cri.smooth} would
imply smoothness of the teardrop algebra ${\mathcal O}(\mathbb{WP}_q(1,l))$.
This arguing leads~to:
\begin{theorem}
${\mathcal O}(S^{3}_q)$ is a~strongly ${\mathbb Z}_l$-graded algebra with the degree-zero subalgebra
isomorphic to ${\mathcal O}(L_q(l;1,l))$.
Consequently, both ${\mathcal O}(L_q(l;1,l))$ and ${\mathcal O}(\mathbb{WP}_q(1,l))$ are homologically
smooth algebras.
\end{theorem}
\begin{proof}
The case $l=1$ is dealt with in~\cite{Kra:Hoc}, the remaining cases follow from
\begin{lemma}
\label{lemma.l.1}
For all integers $l>1$, there exist elements $\omega(1) \in {\mathcal O}(S^{3}_q)_{l-1}\otimes {\mathcal
O}(S^{3}_q)_{1}$ such that $\mu(\omega(1)) = 1$.
\end{lemma}
\begin{proof}
Set:
\begin{gather*}
\omega(1)=x_1\alpha^{l-1}\otimes\alpha^{*l-1}+\sum_{p=1}^{l-1}y_p a^{p-1}\alpha^*\otimes\alpha,
\end{gather*}
where $x_1, y_1, \ldots, y_{l-1} \in {\mathbb C}$ are to be determined and $a=\beta\beta^* = dd^*$.
Then $\omega(1) \in {\mathcal O}(S^{3}_q)_{l-1}\otimes {\mathcal O}(S^{3}_q)_{1}$.
Using~\eqref{su} one f\/inds that
\begin{gather}
\label{aa*}
\alpha^m\alpha^{*m}=\Poch{a}{q^{2}}m=:\sum_{p=0}^m c^m_p a^p,
\end{gather}
where $c^m_p$ are the appropriate $q$-binomial coef\/f\/icients (def\/ined by the second equality
in~\eqref{aa*}).
In view of~\eqref{su}, the condition $\mu(\omega(1)) =1$ leads to
\begin{gather*}
x_1\sum_{p=0}^{l-1}c^{l-1}_p a^p+\sum_{p=0}^{l-2}y_{p+1}a^{p}-q^{-2}\sum_{p=1}^{l-1}y_p a^{p}=1.
\end{gather*}
By comparing the powers of $a$, this is converted into an inhomogeneous system of $l$ equations with
unknown $x_1, y_1, \ldots, y_{l-1}$, whose determinant is
\begin{gather*}
\Delta_{l:1}=\left|
\begin{array}{c c c c c c c}
c^{l-1}_0& 1& 0& 0& \cdots& 0& 0
\\[6pt]
c^{l-1}_1& -q^{-2}& 1& 0& \cdots& 0& 0
\\[6pt]
c^{l-1}_2& 0& -q^{-2}& 1& \cdots& 0& 0
\\[2pt]
\vdots& & & \vdots
\\[2pt]
c^{l-1}_{l-2}& 0& 0& 0& \cdots& -q^{-2}& 1
\\[6pt]
c^{l-1}_2& 0& 0& 0& \cdots& 0& -q^{-2}
\end{array}
\right|.
\end{gather*}
{}$\Delta_{l:1}$ can be evaluated by expanding by the f\/irst column to give
\begin{gather*}
\Delta_{l:1}=(-1)^{l-1}\big(q^{-2(l-1)}c^{l-1}_0+q^{-2(l-2)}c^{l-1}_1+\cdots+c^{l-1}_{l-1}\big)
=\big({-}q^2\big)^{l-1}\prod_{p=1}^{l-1}\big(1-q^{2p}\big)\neq0.
\end{gather*}
The f\/inal equality follows from the def\/inition of the $q$-binomial coef\/f\/icients~\eqref{aa*}.
This proves the existence of $\omega(1)$ as stated.
\end{proof}

Since 1 is a~generator of ${\mathbb Z}_l$, Lemma~\ref{lemma.l.1} ensures the existence of mappings $\omega:
{\mathbb Z}_l\to {\mathcal O}(S^{3}_q)\otimes {\mathcal O}(S^{3}_q)$ that satisfy conditions $(a)$ and $(b)$ in
Lemma~\ref{lem.strong}.
Hence ${\mathcal O}(S^{3}_q)$ is a~strongly ${\mathbb Z}_l$-graded algebra, and the second assertion of the
theorem follows by Criterion~\ref{cri.smooth}.
\end{proof}

Therefore, for any $q\in (0,1)$ the action of ${\mathbb Z}_l$ on the quantum three-sphere described above
is free despite the fact that the corresponding action on the classical level is not free (unless,
obviously, $l=1$).
The f\/ixed points correspond to a~manifold rather than an orbifold.

\subsection[Odd weighted real projective planes ${\mathbb R} {\mathbb P}_q^2(l;-)$ and quantum Seifert lens
spaces]{Odd weighted real projective planes $\boldsymbol{{\mathbb R} {\mathbb P}_q^2(l;-)}$\\ and quantum Seifert lens spaces}

For a~positive integer $l$, the coordinate $*$-algebra ${\mathcal O}({\mathbb R}
{\mathbb P}_q^2(l;-))$ of the odd quantum weighted real projective plane is generated by $a$, $b$, $c_-$
which satisfy the relations:
\begin{gather*}
a=a^*,
\qquad
ab=q^{-2l}ba,
\qquad
ac_{-}=q^{-4l}c_-a,
\qquad
b^2=q^{3l}ac_-,
\qquad
bc_-=q^{-2l}c_-b,
\\
bb^*=q^{2l}a\Poch{a}{q^{2}}l,
\qquad
b^*b=a\Poch{q^{-2}a}{q^{-2}}l,
\qquad
b^*c_-=q^{-l}\Poch{q^{-2}a}{q^{-2}}l b,
\\
c_-b^*=q^{l}b\Poch{a}{q^{2}}l,
\qquad
c_-c_-^*=\Poch{a}{q^{2}}{2l},
\qquad
c_-^*c_-=\Poch{q^{-2}a}{q^{-2}}{2l},
\end{gather*}
see~\cite{Brz:Sei}.
To prove homological smoothness of ${\mathcal O}({\mathbb R} {\mathbb P}_q^2(l;-))$ we make use of
Criterion~\ref{cri.smooth} and build a~tower of strongly graded algebras with ${\mathcal O}({\mathbb R}
{\mathbb P}_q^2(l;-))$ as the foundation.

The coordinate $*$-algebra of the quantum Seifert manifold ${\mathcal O}(\Sigma_q^3)$ is generated by
a~central unitary $\xi$ and elements $\zeta_0$, $\zeta_1$ such that
\begin{gather}
\label{sigma.rel}
\zeta_0\zeta_1=q\zeta_1\zeta_0,
\qquad
\zeta_0\zeta_0^*=\zeta_0^*\zeta_0+\big(q^{-2}-1\big)\zeta^2_1\xi,
\qquad
\zeta_0\zeta_0^*+\zeta_1^2\xi=1,
\qquad
\zeta_1^*=\zeta_1\xi.
\end{gather}
It has been shown in~\cite[proof of Proposition~5.2]{BrzZie:pri} that ${\mathcal O}(\Sigma_q^3)$ can be
understood as the degree-zero part of a~${\mathbb Z}_2$-grading of ${\mathcal O}(S_q^2)[u,u^{-1}]$, where
$u^{-1}=u^*$ and ${\mathcal O}(S_q^2)$ is the coordinate $*$-algebra of the equatorial Podle\'s
sphere~\cite{Pod:sph}, generated by $z_0$ and self-adjoint $z_1$ such that
\begin{gather}
\label{sphere}
z_0z_1=q z_1z_0,
\qquad
z_0z_0^*=z_0^*z_0+\big(q^{-2}-1\big)z_1^2,
\qquad
z_0z_0^*+z_1^2=1.
\end{gather}
The ${\mathbb Z}_2$-grading of ${\mathcal O}(S_q^2)[u,u^{-1}]$ is determined by setting, for all monomials
$w$ of degree $k$ in the basis $\{z_0^rz_1^s, z_0^{*r}z_1^s \; | \: r,s\in {\mathbb N}\}$ of ${\mathcal   
O}(S_q^2)$, $\deg(wu^m) = (k+m)\mod 2$.
{}${\mathcal O}(\Sigma_q^3)$ can be identif\/ied with the degree-zero part of ${\mathcal
O}(S_q^2)[u,u^{-1}]$ by $*$-embedding $\zeta_i\mapsto z_iu$, $\xi\mapsto u^{-2}$.
Thanks to the last of equations~\eqref{sphere}, the function
\begin{gather*}
\omega:{\mathbb Z}_2\to{\mathcal O}\big(S_q^2\big)\big[u,u^{-1}\big]\otimes{\mathcal O}\big(S_q^2\big)\big[u,u^{-1}\big],
\qquad
0\mapsto1\otimes1,
\qquad
1\mapsto z_0\otimes z_0^*+z_1\otimes z_1,
\end{gather*}
satisf\/ies conditions (a) and (b) in Lemma~\ref{lem.strong}, hence ${\mathcal O}(S_q^2)[u,u^{-1}]$ is
a~strongly ${\mathbb Z}_2$-graded algebra.
Since ${\mathcal E}{\mathcal O}(S_q^3)$ is Noetherian, and there is a~surjective $*$-algebra homomorphism
${\mathcal O}(S_q^3)\to {\mathcal O}(S_q^2)$, $\alpha\mapsto z_0$, $\beta\mapsto z_1^*$, both ${\mathcal
E}{\mathcal O}(S_q^2)$ and ${\mathcal E}{\mathcal O}(S_q^2)[u,u^{-1}]$ and hence also ${\mathcal
E}{\mathcal O}(\Sigma_q^3)$ are Noetherian.

As explained in~\cite{BrzFai:rps}, ${\mathcal O}(\Sigma_q^3)$ is a~${\mathbb Z}_{l}$-graded algebra with
grading given by
\begin{gather*}
\deg(\zeta_0)=1,
\qquad
\deg(\zeta_0^*)=l-1,
\qquad
\deg(\zeta_1)=\deg(\xi)=0.
\end{gather*}
The degree-zero part of ${\mathcal O}(\Sigma^3_q)$ is isomorphic to the $*$-algebra ${\mathcal
O}(\Sigma^3_q(l;-))$ generated by $x$, $y$ and central unitary $z$ subject to the following relations
\begin{gather*}
y^*=yz,
\qquad
xy=q^l yx,
\qquad
xx^*=\Poch{y^2z}{q^{2}}l,
\qquad
x^*x=\Poch{q^{-2}y^2z}{q^{-2}}l.
\end{gather*}
The embedding of ${\mathcal O}(\Sigma^3_q(l;-))$ into ${\mathcal O}(\Sigma_q^3)$ is given by $x\mapsto
\zeta_0^l$, $y\mapsto \zeta_1$ and $z\mapsto \xi$.
The similarity of relations~\eqref{sigma.rel} and~\eqref{su} leads immediately to equations~\eqref{aa*}
with $\alpha$ replaced by $\zeta_0$ and $a= \zeta_1^2\xi$.
This allows one to use the same arguments as in Lemma~\ref{lemma.l.1} to prove that there exist $x_1, y_1,
\ldots, y_{l-1}\in {\mathbb C}$ such that
\begin{gather*}
\omega(1)=x_1\zeta_0^{l-1}\otimes\zeta_0^{*l-1}+\sum_{i=1}^{l-1}y_ia^{i-1}\zeta_0^{*}\otimes\zeta_0^{}
\in{\mathcal O}\big(\Sigma_q^3\big)_{l-1}\otimes{\mathcal O}\big(\Sigma_q^3\big)_l,
\end{gather*}
has the required property $\mu(\omega(1)) =1$.
Therefore, ${\mathcal O}(\Sigma^3_q)$ is a~strongly graded ${\mathbb Z}_l$-algebra.

Finally, it is proven in~\cite{BrzFai:rps} that ${\mathcal O}(\Sigma^3_q(l;-))$ is a~strongly ${\mathbb
Z}$-graded algebra with grading given by $\deg(x)=\deg(y)=1$, $\deg(x^*) =-1$ and $\deg(z) = -2$.
The degree-zero subalgebra of ${\mathcal O}(\Sigma^3_q(l;-))$ can be identif\/ied with the coordinate
algebra of weighted real projective plane ${\mathcal O}({\mathbb R}{\mathbb P}^2_q(l;-))$ via the map
$a\mapsto y^2z$, $b\mapsto xyz$ and $c_-\mapsto x^2z$.

Summarizing, we have presented in this section a~tower of $*$-algebras
\begin{gather}
\label{seq}
{\mathcal O}\big({\mathbb R}{\mathbb P}^2_q(l;-)\big)\hookrightarrow{\mathcal O}
\big(\Sigma^3_q(l;-)\big)\hookrightarrow{\mathcal O}\big(\Sigma^3_q\big)\hookrightarrow{\mathcal O}\big(S_q^2\big)\big[u,u^{-1}\big].
\end{gather}
The second, third and fourth terms are strongly group graded algebras.
Each antecedent term is the degree-zero part of the subsequent one.
Since the enveloping algebra of ${\mathcal O}(S_q^2)[u,u^{-1}]$ is Noetherian, so are the enveloping
algebras of all its predecessors.
By~\cite[Corollary~4.6]{Liu:hom} the global dimension of ${\mathcal E}{\mathcal O}(S_q^2)$ is f\/inite,
hence so is the global dimension of ${\mathcal E}{\mathcal O}(S_q^2)[u,u^{-1}]$, and, by
Criterion~\ref{cri.smooth}, the global dimensions of enveloping algebras of all its predecessors
in~\eqref{seq}.
This proves the following
\begin{theorem}
The algebras ${\mathcal O}(\Sigma^3_q)$, ${\mathcal O}({\mathbb R}{\mathbb P}^2_q(l;-))$
and ${\mathcal O}(\Sigma^3_q(l;-)) $ are homologically smooth.
\end{theorem}

\subsection[Quantum real weighted projective planes ${\mathbb R}{\mathbb P}_q^2(l;+)$ and teardrops
(revisited)]{Quantum real weighted projective planes $\boldsymbol{{\mathbb R}{\mathbb P}_q^2(l;+)}$\\ and teardrops (revisited)}

Let $k$ be a~natural number and $l$ be a~positive integer.
Write ${\mathcal A}(k,l)$ for the $*$-algebra generated by $a$ and $b$ subject to the following relations
\begin{gather*}
a^*=a,
\qquad
ab=q^{-2kl}ba,
\qquad
bb^*=q^{2kl}a^k\Poch{a}{q^{2}}l,
\qquad
b^*b=a^k\Poch{q^{-2}a}{q^{-2}}l.
\end{gather*}
If $k$ and $l$ are coprime then ${\mathcal A}(k,l)$ is the coordinate algebra of the quantum weighted
projective line or the {\em quantum spindle} ${\mathcal O}(\mathbb{WP}_q(k,l))$ introduced
in~\cite{BrzFai:tea}.
The special case $k=1$ is simply the {\em quantum teardrop}; see Section~\ref{sec.lens}.
For $l$ odd, ${\mathcal A}(0,l)$ is the coordinate algebra of the quantum weighted even real projective
plane ${\mathcal O}({\mathbb R}{\mathbb P}_q^2(l;+))$ introduced in~\cite{Brz:Sei}.
The following theorem is a~consequence of Criterion~\ref{cri.Weyl}.
\begin{theorem}
The algebras ${\mathcal A}(k,l)$ are homologically smooth $($of dimension $2)$ if and only if $k=0,1$.
\end{theorem}
\begin{proof}
We only need to observe that each ${\mathcal A}(k,l)$ is a~generalized Weyl algebra over the polynomial
algebra ${\mathbb C}[a]$ given by the automorphism $\pi(a) = q^{2l} a$, element $p=a^k
\prod\limits_{m=1}^l(1-q^{-2m}a)$ and generators $x_-=b$, $x_+ =b^*$.
Since $p$ has no multiple roots if and only if $k=0,1$, the assertion follows by Criterion~\ref{cri.Weyl}.
\end{proof}

Furthermore, for $k=0,1$, ${\mathcal A}(k,l)$ are twisted Calabi--Yau algebras with the twisting
automorphism $\nu(b) = q^{-2l}b$, $\nu(b^*)= q^{2l}b^*$ and $\nu(a)=a$.
Hence they enjoy the Poincar\'e duality in the sense of Van den Bergh~\cite{Van:rel}.

\subsection*{Acknowledgements}

I would like to thank Ulrich Kr\"ahmer for discussions, Li-Yu Liu for
bringing reference~\cite{Liu:hom} to my attention, and Piotr~M.~Hajac and the referees for helpful comments.
I am grateful to Fields Institute for Research in Mathematical Sciences in Toronto, where these results were f\/irst presented, for creating excellent research environment and for support.

\pdfbookmark[1]{References}{ref}
\LastPageEnding

\end{document}